\documentclass[reqno,12pt,draft]{amsart}
\usepackage{amsfonts,amsmath,amssymb,amscd}
\usepackage[russian]{babel}
\input epsf
\long\def\comment#1\endcomment{}

\begin{document}

\centerline{\uppercase{\bf Размышления о признании геометрии Лобачевского}
\footnote{
Благодарим  И. Изместьева,
А. Петрунина
и членов редколлегии журнала `Математическое Просвещение' за полезные обсуждения.
Особо благодарим А. Сосинского за разрешение использовать часть написанного им текста.}
}
\bigskip
 \centerline{\bf В.В. Прасолов и А.Б. Скопенков
\footnote{Инфо: www.mccme.ru/prasolov, www.mccme.ru/\~{ }skopenko.
А. Скопенков частично поддержан грантом фонда Саймонса.}
}




\bigskip
Признание неевклидовой геометрии происходило непросто.
Часто споры вокруг нее представляют упрощенно, в черно-белом цвете:
`невежи и консерваторы осмеивали гениев и новаторов'.
И хотя в этом имеется доля истины, дело здесь обстоит сложнее, и до поры
до времени для   неприятия этой новой теории были довольно веские основания.
Мы попытаемся преодолеть этот стереотип, в частности, объяснить, {\it почему неевклидова геометрия  была принята не сразу}, и показать,
что для ее признания было важно не только открытие ее моделей (\S1), но и {\it появление ее приложений к другим областям математики} (\S2).

Приводимые размышления о приложениях напрямую касаются \emph{современной} математики и связаны
с важнейшими практическими вопросами,
например: как математику выбирать направления исследования.
Эти мысли знакомы
некоторым профессиональным математикам.
Но мы надеемся, что они будут новы и интересны и для тех, кто изучает математику (на разных уровнях),
и для тех, кто интересуется математикой как частью культуры.

Мы не рассматриваем {\it предысторию} геометрии Лобачевского (т.е. труды Аристотеля, Омара Хайяма,
Швейкарта, Тауринуса, Ламберта, Саккери и других).
Мы лишь коротко пишем об истории открытия геометрии Лобачевского и ее моделей
(\S1; если вам знакома эта история, то вы можете пропустить \S1).
Подробнее см. [Pa], [Gr] и [W].
Основы самой геометрии Лобачевского можно изучить, например, по [Pr, Pe].

\section{Проблема Пятого постулата и геометрия Лобачевского}


Пятый постулат Евклида (аксиома параллельных) в более простой эквивалентной форме гласит, что

{\it через точку вне данной прямой проходит не более одной прямой, не пересекающей  данную прямую.}

Эта аксиома формулируется сложнее, чем другие аксиомы (ибо она требует рассмотрения бесконечной прямой).
Поэтому еще Аристотель рассуждал о возможности вывести ее из остальных аксиом.
(Точнее, Аристотель рассуждал об эквивалентном утверждении о сумме углов треугольника.)
Чтобы сделать это, Карл-Фридрих Гаусс, Янош Бойяи и Николай Иванович Лобачевский, как и их предшественники,
рассуждали от противного.
Они предположили, что через точку вне данной прямой проходит более, чем  одна прямая, не пересекающая данную, и стремились отсюда получить противоречие.
Но следствия из этой гипотезы, вместо того, чтобы привести к противоречию, постепенно выстраивались в очень стройную и богатую, хотя и крайне необычную, теорию.
Появилось подозрение, что эта теория так же логически безупречна, как и евклидова, а потом и уверенность в этом.
Такая уверенность пришла к Гауссу, Бойяи и Лобачевскому в 1810х-1820х годах.
Первая публикация принадлежит Лобачевскому.

Однако в начале 19-го века не было ясно, как, развивая неевклидову геометрию, можно решить двухтысячелетнюю
проблему о доказуемости Пятого постулата Евклида.


Проблема Пятого постулата была решена после появления {\it моделей} неевклидовой геометрии, придумать которые помогли труды Лобачевского.
Поясним читателю важное понятие модели.
Предположим, что дана аксиоматическая теория, непротиворечивость которой мы хотим установить.
(Например, неевклидова геометрия.)
В этой теории, кроме аксиом, должны быть \emph{неопределяемые понятия} --- это понятия,
определения которых в теории не даются.
Такие понятия в строгой аксиоматической теории необходимы, иначе в определениях будет порочный круг (так же как необходимы недоказываемые \emph{аксиомы}, без которых в доказательствах будет порочный круг).

Попробуем  обратиться к другому разделу математики, в непротиворечивости которого мы не сомневаемся
(например, к евклидовой геометрии).
И построим в нем \emph{модель} данной теории, т.е. переведем неопределяемые понятия данной теории в конкретные
термины выбранного раздела математики так, чтобы аксиомы теории (верней, их перевод) превратились в истинные высказывания этого раздела.
Тогда, если имеется противоречие в теории (например, в геометрии Лобачевского), то соответствующее противоречие
можно получить и в том разделе математики, где реализована наша модель (например, в евклидовой геометрии).
Если же в модели противоречия нет, то его нет и в исходной теории.
Установленную таким способом непротиворечивость называют {\it относительной} потому, что она может быть получена лишь в предположении, что тот раздел математики, внутри которого построена  модель, сам является непротиворечивым.
Так что из существования указанной модели следует, что если есть противоречие в геометрии Лобачевского, то оно есть и в геометрии Евклида.
\footnote{Это не дает ответа на важнейший исходный вопрос: есть ли противоречие в геометрии Лобачевского?
Но, с одной стороны, в непротиворечивости геометрии Евклида все уверены.
С другой стороны, появление моделей и доказательство относительной непротиворечивости подняло вопрос о
формальном доказательстве непротиворечивости геометрии Евклида (и других `классических' теорий).
Однако даже формализация этого вопроса нетривиальна, и останется проблема корректности самих методов доказательства.
Мы благодарны Д. Мусатову за обсуждение этой важной темы и надеемся, что она будет освещена в популярной литературе.}

В своих работах Лобачевский писал, что у обитателей пространства с его геометрией не было бы проблем с евклидовой геометрией, поскольку геометрия на {\it орисфере} евклидова.
({\it Орисферой} в геометрии Лобачевского называется `предельная  фигура', полученная из сферы устремлением радиуса к бесконечности. При этом предполагается, что сфера проходит через некоторую фиксированную точку, а ее центр уходит
в бесконечность по фиксированному лучу, выходящему из этой точки.
В отличие от евклидовой геометрии, эта фигура не  совпадает с плоскостью.)
Видимо, он интуитивно пришел к идее построения модели, но не сформулировал этого явно.
Главной трудностью в доказательстве непротиворечивости неевклидовой геометрии в те времена было
отсутствие явно высказанной идеи модели.
А в наше время идея построения модели для доказательства (относительной) непротиворечивости общеизвестна и даже банальна.


Первая модель неевклидовой геометрии была построена итальянским математиком Эудженио Бельтрами [K2].
Тем самым было доказано, что Пятый постулат невозможно вывести из остальных аксиом.

Сначала, в 1868 г., Бельтрами построил модель {\it малой части} плоскости Лобачевского [B].
Приведем идею его построения.
Трактриса (кривая погони) --- плоская кривая, для которой длина отрезка касательной от точки касания до точки пересечения с фиксированной прямой является постоянной величиной.
Такую линию описывает предмет, волочащийся на веревке заданной длины за точкой, движущейся по оси абсцисс.
Поверхность вращения трактрисы называется {\it поверхностью Бельтрами} (хотя она была открыта до Бельтрами).
\footnote{Тригонометрические формулы для поверхности Бельтрами (и других поверхностей постоянной отрицательной кривизны) были получены Фердинандом Миндингом (прибалтийским учеником Гаусса) [M] в 1840 г.
Они совершенно идентичны формулам геометрии Лобачевского, открытым Лобачевским!}
{\it Геодезической} на поверхности Бельтрами (и на произвольной поверхности) называется любая линия, достаточно малые   дуги которой являются на этой поверхности кратчайшими путями между их концами.
Бельтрами доказал, что геометрия на {\it малой части} поверхности Бельтрами такая же, как на {\it малой части} плоскости Лобачевского.
Т.е. что перевод
\begin{eqnarray*}
\text{часть плоскости Лобачевского}\to\text{ поверхность Бельтрами}\\
\text{точка  }\to\text{ точка этой поверхности}\\
\text{прямая }\to\text{ геодезическая на этой поверхности}
\end{eqnarray*}
доставляет модель части плоскости Лобачевского в виде поверхности Бельтрами.

Затем Бельтрами построил модель {\it всей} плоскости Лобачевского (весьма неожиданную и простую!) [K, Pr].
Важные свойства этой модели, на которых сейчас основаны ее приложения к другим областям математики, были обнаружены Феликсом Клейном в 1871
(связь с проективной геометрией, а именно --- с проективными преобразованиями, сохраняющими круг).
Клейн использовал важную идею Артура Кэли (1859) о связи расстояний на сфере и {\it двойных отношений} из
{\it проективной геометрии}.
Поэтому модель, построенная Бельтрами, называется моделью Кэли-Клейна или Бельтрами-Клейна.

Некоторые другие модели (1883) называют именем Анри Пуанкаре.
Хотя Бельтрами писал о них, именно Пуанкаре указал на связь с другими областями математики
(см. подробнее \S2).

Интересно, что когда велись споры вокруг геометрии Лобачевского, другая неевклидова геометрия уже была давно общепризнана.
Это геометрия звездного неба или поверхности земного шара --- сферическая геометрия.
Ввиду реальности изучаемого объекта и его огромного значения для астрономии и мореплавания
никаких насмешек и борьбы за признание сферической геометрии просто не было.
Однако сферическая геометрия до Клейна (может быть --- до Римана) воспринималась не как отдельная геометрия, а как
часть трехмерной евклидовой геометрии.
\footnote{Пятый постулат пытались вывести из остальных аксиом Евклида.
На сфере не все они выполнены.
Поэтому сферическая геометрия не воспринималась как неевклидова и как имеющая тесную связь с геометрией Лобачевского.}
Поэтому не было необходимости в построении модели, как не было и самого понятия модели.



\section{Признание геометрии Лобачевского и приложения}

{\bf Почему Гаусс не публиковал свои работы по неевклидовой геометрии?}

Гаусс первый осознал, что неевклидова геометрия имеет право на существование.
Почему же он не продолжил свои занятия этой областью?
Часто приходится читать, что Гаусс не публиковал из `интеллектуальной трусости'  ---
он боялся, что его поднимут на смех
\footnote{Кстати, именно это случилось с Лобачевским, о котором с насмешкой писал Чернышевский [Gi, с. 376].
Да и в журнале `Сын отечества' (известном еще травлей Пушкина) труды Лобачевского были грубо осмеяны.}.
Это подтверждается его письмом 1829 года Бесселю, где Гаусс признает состоятельность неевклидовой геометрии, но подчеркивает, что не объявляет это публично во избежание криков `беотийцев' (т. е. невежественных людей).
В какой-то мере, возможно, это было и так.
Но нам представляется, что более существенны незавершенность попытки решить проблему Пятого постулата, а также
следующие весьма достойные причины.

В начале своего научного пути Гаусс занимался `чистой' математикой (в частности, теорией чисел, замечательные приложения которой были обнаружены лишь позднее).
Начиная с 1816 года, когда он обосновался в Геттингене, он занимался разделами математики (теорией вероятностей, векторным анализом, дифференциальной геометрией и др.), связанными с практическими приложениями (астрономией, геодезией, магнетизмом).
Значит, Гаусс считал более важными для себя области математики, связанные с изучением реального мира.
Это косвенно подтверждает и отрывок из его рецензии [Ga]
(на неудачную попытку доказательства пятого постулата; перевод авторов):


{\it Большая часть [рецензируемой] работы касается утверждения, что, вопреки Канту,
достоверность геометрии базируется не на лицезрении, а на
определениях и логических правилах вывода.
Кант вовсе не хотел отрицать, что эти логические вспомогательные средства все
более и более используются для описания геометрических истин и связей между ними.
Однако любой человек, знакомый с сущностью геометрии, согласится, что
логические средства сами по себе не позволяют ничего получить, а дают лишь
пустоцвет, если всюду не властвует плодотворное живое лицезрение предметов.}


Лобачевский, в отличие от Гаусса, не имеет ярких результатов, связанных с приложениями.
Однако его деятельность на посту ректора Казанского университета не менее важна для реальной жизни, чем
прикладные исследования.
При этом вызывает уважение то, что с одной стороны, сам он продолжал заниматься неевклидовой геометрией,
а с другой стороны, не использовал своего высокого служебного положения для продвижения своих исследований.
Будучи ректором Казанского университета, он мог бы основать научную школу по неевклидовой геометрии и издавать собственный журнал, не дожидаясь ее международного признания.
В этом проявилось отличие большого ученого и порядочного человека от профана или карьериста,
пытающегося любой ценой продвинуть свои идеи.


\smallskip
{\bf Астрономические  наблюдения Лобачевского }

Как указано в [V], Лобачевский производил астрономические наблюдения с целью проверить, равна ли сумма углов треугольника $180^{\circ}$ (что эквивалентно пятому постулату Евклида) или меньше, как в `воображаемой геометрии'.
\footnote{Как указано в [K2, с. 26-27], Гаусс производил аналогичные наблюдения на поверхности Земли.
Гаусс проводил картографические исследования и, возможно, при этом измерял большие треугольники
(образованные вершинами гор) для соединения измерений, выполненных для различных участков карт.
Возможно, это и привело к появлению легенды о том, что он проводил этот дорогостоящий эксперимент для удовлетворения своего любопытства, связанного с геометрией Лобачевского.}

Здесь мы подходим к одному из ключевых вопросов философии науки: что такое геометрия, о чем эта наука?
Современники и предшественники Лобачевского (да и он сам до поры до времени) считали, что трехмерная
геометрия Евклида --- это учение о физическом пространстве нашего мира.
Но у Лобачевского в какой-то момент забрезжила мысль --- а может быть, его геометрия вовсе не такая уж
`воображаемая', и именно она, а не геометрия Евклида, и определяет структуру нашего пространства?

Астрономические измерения Лобачевского (мы их здесь не описываем, но заинтересованный читатель может найти рассказ о них в книге [V]) не привели  к ответу: сумма углов получилась меньше $180^{\circ}$, но отличие от $180^{\circ}$ не превысило ошибку (точность) измерений.
Вопрос о том, какая из геометрий и есть геометрия нашей Вселенной, так и остался висеть в воздухе.

\smallskip
{\bf Первые приложения геометрии Лобачевского}

Важнейшую роль в признании геометрии Лобачевского сыграло открытие моделей (\S1).
Были также важны и логический анализ оснований геометрии (Паш, Гильберт и др.), вызванный появлением геометрии Лобачевского и ее моделей, и `эрлангенская программа' Клейна [K1] вместе с теорией Ли {\it непрерывных групп}.

Большое значение имело также открытие приложений к другим областям математики.

Сам Лобачевский вычислил множество определенных интегралов, интерпретируя их как объемы различных тел в неевклидовом пространстве.
Однако эти применения были раскритикованы М.В. Остроградским.
Он указывал, что один из двух интегралов, вычисленный Лобачевским, известен,
а второй неверен [Pa, стр. 14 внизу].
Ввиду критики Остроградского это приложение геометрии Лобачевского не сильно способствовало ее признанию.
(Мы не обсуждаем здесь вопрос о том, в какой степени эта критика была справедливой.)

Одним из первых приложений геометрии Лобачевского к другим областям математики была теория \emph{автоморфных функций},
разработанная Пуанкаре в 1881-84 годы [K2].
В простейшем случае это функции комплексного переменного, определённые в верхней полуплоскости и инвариантные относительно некоторого множества дробно-линейных преобразований
вида $z\mapsto\dfrac{az+b}{cz+d}$ с вещественными коэффициентами $a$, $b$, $c$ и $d$.
Пуанкаре сначала выясняет, как устроены {\it фундаментальные области}
\footnote{Читатель, не знающий, что это такое, может пропустить этот абзац без ущерба для понимания дальнейшего.}
таких преобразований, а затем строит сами функции с помощью рядов.
Пуанкаре обнаруживает, что фундаментальные области заполняют верхнюю полуплоскость, причeм их размеры уменьшаются при
приближении к границе --- вещественной прямой.
Это напомнило ему геометрию Лобачевского, и неожиданно пришла идея, что эти преобразования совпадают с движениями неевклидовой
геометрии. Правильность этой идеи Пуанкаре вскоре легко проверил.

Пуанкаре отмечал, что геометрия Лобачевского служила ему в его исследованиях руководящей нитью, но он избегал использовать
ее в своем изложении, потому что она была в то время мало знакома математикам.
Этому знакомству весьма поспособствовали исследования Пуанкаре, показавшие, что геометрия Лобачевского
может иметь приложения и во вполне классических областях математики.

Некоторую роль в признании геометрии Лобачевского могло сыграть также открытие в начале 20 века ее приложений
к физике, точнее, к теории относительности [DSS].
Оказалось, что {\it пространство скоростей} в этой теории имеет геометрию Лобачевского
(иными словами, `совпадает' с моделью Бельтрами-Клейна).

\smallskip
{\bf Заключение: `естественнонаучный' и `философский' аспекты математики}

Лобачевский, Гаусс и Бойяи высказали две важные идеи.
Во-первых, логически мыслима не только евклидова геометрия, но и другие геометрии.
Во-вторых, эти другие геометрии в принципе могут отражать строение реального пространства.
К сожалению, в то время эти две разные идеи не были четко отделены друг от друга.
Споры вокруг неевклидовой геометрии помогли математикам четко выделить два разных аспекта своей науки.
Первый --- изучение систем аксиом и, вообще, формальных конструкций; он ближе к философии.
Второй --- математическое изучение реального мира; он ближе к естественным наукам (в первую очередь, к физике).
Эти аспекты взаимосвязаны: математическое изучение реального мира порождает математические реальности,
уже не так непосредственно связанные с реальным миром.
Нам близки естественнонаучные позиции Гаусса, Пуанкаре, Колмогорова и Арнольда: математическое изучение реального мира --- важнейший аспект математики, но изучение формальных конструкций также необходимо.

Конечно, эти идеи важны и достойны более детального обсуждения.
Однако оно не входит в нашу цель, мы хотели лишь еще раз обозначить их.
Будем рады, если это приведет к последующим обсуждениям и публикациям.

\small

\section*{Дополнение: о публикациях}

Почему Гаусс не опубликовал свои размышления, написано в начале \S2.
В отличие от Гаусса, Я. Бойяи был готов публиковать свою работу сразу.
Но это оказалось непросто: где найти издателя, готового опубликовать столь необычный труд?
На выручку пришел его отец,      включив эту работу в свою книгу по геометрии в виде приложения.
В 1832 году двухтомная книга Фаркаша Бойяи, \emph{Tentamen}, содержащая знаменитый сегодня \emph{Аппендикс} Яноша, выходит в свет.


Лобачевскому тоже непросто было публиковать свои работы в изданиях, читаемых большим количеством математиков.
Но, в отличие от Гаусса и Бойяи, он сумел это сделать достаточно полным образом.
Проблема была в другом: на публикации Лобачевского по-русски никто из серьезных ученых не обратил внимания
(например, даже выдающийся русский математик Буняковский в своей работе о теории параллельных вовсе не упомянул Лобачевского), да и никто из французских математиков не обратил внимание на его последний труд  (\emph{Pang\'eometrie}).
Не было и немецких читателей у его \emph{Geometrische Untersuchungen}.
Кроме одного, но зато какого --- эту небольшую книгу прочитал Гаусс и был потрясен.
Известно, что он стал учить русский язык, возможно, чтобы прочитать более ранние публикации Лобачевского в Казанском журнале.
Гаусс добился избрания Лобачевского членом-корреспондентом Г\"eттингенского королевского научного общества.


{\it Cовременному} математику (или физику) гораздо проще сделать свои исследования доступными мировому сообществу ученых (и, тем самым, защитить свой приоритет).
Имеются научные конференции и конгрессы.
Кроме того, любой автор научной работы может
выложить ее на международный сервер http://arxiv.org (архив).
Получить для этого рекомендацию кого-то, кто уже имеет выложенные статьи, нетрудно.
Выкладывание работы в архив накладывает определенную ответственность на ее автора:
его репутация может пострадать, если работа ошибочная или недостаточно обоснованная.
Но все же из-за свободы выкладывания в архив в нем появляется много `мусора'.
Поэтому само по себе выкладывание в архив не гарантирует, что работу прочитают (и притом не считается официальной публикацией).
В первом приближении, больше шансов быть прочитанными имеют работы,

$\bullet$ ясно написанные;

$\bullet$ представленные на конференциях и семинарах;

$\bullet$ известных авторов;

$\bullet$ по модной тематике;

$\bullet$ не только выложенные в архив, но и параллельно опубликованные в хороших журналах.

Несмотря на все `но', в наше время есть гораздо больше возможностей для распространения новых научных идей,
чем во времена Гаусса, Бойяи и Лобачевского.

\normalsize

\bigskip
{\bf БИБЛИОГРАФИЯ}


[B] Э.Бельтрами. Опыт интерпретации неевклидовой геометрии. В сб.: `Об основаниях геометрии.'
М.: Гостехиздат. 1956. С. 180 - 212.

[DSS] В.Н. Дубровский, Я.А. Смородинский, Е.Л. Сурков, Релятивистский мир, М., Наука, 1984.
http://ilib.mccme.ru/djvu/bib-kvant/kvant34.htm

[Ga] C.F.Gauss. Werke, Bd 4, Goetingen, 1873. S. 364-368.
(aus G\"ottingsche gelehrte Anzeigen).

[Gi] Гиндикин, Рассказы о физиках и математиках, М.: МЦНМО, 2001.

[Gr] Gray J. Worlds out of nothing. A course in the history of geometry in the
19th century (Springer, 2007)(381s)
http://gen.lib.rus.ec, зеркало: http://free-books.dontexist.com


[K1] Ф.Клейн. Сравнительное обозрение новейших геометрических исследований (`Эрлангенская программа'). В сб.: `Об основаниях геометрии.' М.: Гостехиздат. 1956. С. 399-434.

[K2] Ф. Клейн, Лекции о развитии математики в XIX столетии, М: Наука, 1989.



[M] F. Minding, Beitr\"age zur Theorie der k\"urzesten Linien auf krummen Fl\"achen,
Journal f\"ur die reine und angewandte Mathematik (Crelle's Journal), 20 (1840), 323-327.
\linebreak
http://www.deepdyve.com/lp/de-gruyter/

beitr-ge-zur-theorie-der-k-rzesten-linien-auf-krummen-fl-chen-060ar0lxeG

Рус. перевод:  Ф.Миндинг. Дополнения к теории кратчайших линий на кривых поверхностях.
В сб.: `Об основаниях геометрии.' М.: Гостехиздат. 1956. С. 176-179.



[Pa] А. Пападопулос,
О гиперболической геометрии и истории ее признания,
Матем. Просвещ. 14 (2010), 10-29.

[Pe] A. Petrunin, Euclidean and Hyperbolic Planes; a minimalistic introduction with metric approach.
http://arxiv.org/pdf/1302.1630.pdf

[Pr] В.В.Прасолов. Геометрия Лобачевского (М: МЦНМО, 1995, 2000, 2004),
\linebreak
http://www.mccme.ru/prasolov.

[V] А.В. Васильев, Николай Иванович Лобачевский, М: Наука, 1992.

[W] http://ru.wikipedia.org/wiki/ Геометрия\_Лобачевского

\end{document}